\documentclass[preprint,12pt]{elsarticle}

\usepackage{amssymb}
\usepackage{amsthm}
\usepackage{amsmath}

 \newtheorem{thm}{Theorem}

  \newdefinition{rmk}{Remark} 
  \newproof{pf}{Proof}

\journal{Applied Mathematics Letters}

\begin{document}

\begin{frontmatter}



\title{Optimal properties of tensor product of B-bases\footnote{This work was partially supported through the Spanish research grant PGC2018-096321-B-I00 (MCIU/AEI) and by Gobierno de Arag\'on (E41-17R).}\footnote{Organization: Departamento de Matem\'atica Aplicada/IUMA, Universidad de Zaragoza}}



\author{Jorge Delgado}\ead{jorgedel@unizar.es}
\author{H\'ector Orera} \ead{hectororera@unizar.es}
\author{J. M. Pe\~na} \ead{jmpena@unizar.es}


\begin{abstract}
It is proved the optimal conditioning for the $\infty$-norm of  collocation matrices of the tensor product of normalized B-bases among the tensor product of all normalized totally positive bases of the corresponding space of functions. Bounds for the minimal eigenvalue and singular value and illustrative numerical examples are also included.  
\end{abstract}


\begin{highlights}
\item Minimal conditioning for the $\infty$-norm of collocation matrices of the tensor product
of normalized B-bases is shown.
\item The maximality of the minimal eigenvalue and singular value of collocation matrices 
of the tensor product of normalized B-bases is shown.
\end{highlights}


\end{frontmatter}


\section{Introduction and main results}
\label{}




Given a
system of functions $u = (u_0, \ldots, u_n)$ defined on $I \subseteq
\mathbb{R}$, the {\it collocation matrix} of $u$ at $t_0 < \cdots <
t_m$ in $I$ is given by
$
(u_j(t_i))_{i=0, \ldots, m}^{j=0, \ldots, n}$. If $\sum_{i=0}^n u_i(t) = 1$ for all $t \in
I$, then we say that the system is {\it normalized}. If all collocation matrices of $u$ have all their minors nonnegative, then we say that the system is {\it totally positive} (TP). Normalized totally positive (NTP) systems play a crucial role in Computer Aided Geometric Design because they lead to shape preserving representations. Among all NTP bases of a space, the basis with optimal shape preserving  properties is the {\it normalized B-basis} (\cite{CP,DP}). The Bernstein basis of polynomials and the B-spline basis are examples of normalized B-bases of their corresponding spaces.
In this paper we extend some optimal properties of normalized B-bases given in \cite{DP} to their corresponding tensor products. Recall that, given two systems $u^1=(u_0^1,\ldots,u_m^1)$ and
$u^2=(u_0^2,\ldots,u_n^2)$ of functions defined on
$[a,b]$ and $[c,d]$, respectively, the system $u^1\otimes
u^2:=(u_i^1(x)\cdot
u_j^2(y))_{i=0,\ldots,m}^{j=0,\ldots,n}$ is called a tensor
product system and generates a tensor product surface. 
The \textit{Kronecker product} of two square matrices $A=(a_{ij})_{1\leq i,j\leq m}$ and $B=(b_{ij})_{1\leq i,j\leq n}$, $A\otimes B$, is defined to be the $mn\times mn$ block matrix
\[
A\otimes B=\left(\begin{array}{ccc}
a_{11}B &\cdots   & a_{1m}B\\ 
 \vdots& \ddots &  \vdots\\ 
 a_{m1}B& \cdots &  a_{mm}B
\end{array}\right).
\]
Given the collocation matrices $
B_1:=(u_j^1(x_i))_{0\le i,j\le  m}$ and $
B_2:=(u_j^2(y_i))_{0\le i,j\le  n}$ of $u^1$ and $u^2$, $B_1\otimes B_2$ is the collocation matrix of $u^1\otimes u^2$ at  $((x_i,y_j)_{j=0,\ldots,n})_{i=0,\ldots,m}$. 

Given two square real matrices $A=(a_{ij})_{1\leq i,j \leq n}$ and  $B=(b_{ij})_{1\leq i,j \leq n}$, $A\leq B$ denotes that $a_{ij}\leq b_{ij}$ for all  $i,j$. Given a complex matrix $C=(c_{ij})_{1\leq i,j \leq n}$, $A$ is said to \textit{dominate} $C$ if $|c_{ij}|\leq a_{ij}$ for all $i,j$. 
If matrices $A$ and $B$ are nonsingular, by Corollary 4.2.11 of \cite{topics} we have that $A\otimes B$
is nonsingular and 
\begin{equation}\label{prop.inv}
(A\otimes B)^{-1}=A^{-1}\otimes B^{-1}.
\end{equation}

The next result shows the optimal properties of a collocation matrix of the tensor product of normalized B-bases among all the corresponding collocation matrices of the tensor product of NTP bases of the spaces.

\begin{thm}\label{thm:main}
Let $u^1=(u_0^1,\ldots,u_m^1)$ be an NTP basis on $[a,b]$ of a space of functions $\mathcal{U}_1$, $u^2=(u_0^2,\ldots,u_n^2)$ be an NTP basis on $[c,d]$ of a space of functions $\mathcal{U}_2$ and let $v^1=(v_0^1,\ldots,v_m^1)$ and  $v^2=(v_0^2,\ldots,v_n^2)$ be the normalized B-bases of $\mathcal{U}_1$ and $\mathcal{U}_2$, respectively. Given the increasing sequences of nodes $\textbf{t}=(t_i)_{i=0}^m$ on $[a,b]$ and $\textbf{r}=(r_i)_{i=0}^n$ on $[c,d]$, the nonsingular collocation matrices $A_1$ and $M_1$ of the bases $u^1$ and $v^1$, respectively, at   \textbf{t}, and $A_2$ and $M_2$ 
of the bases $u^2$ and $v^2$, respectively, at  \textbf{r},
 the following properties hold:
\begin{enumerate}[i)]
\item The matrix $|(A_1\otimes A_2)^{-1}|$ dominates $(M_1\otimes M_2)^{-1}$.

\item The minimal eigenvalue (resp., singular value) of $A_1\otimes A_2$ is bounded above by the minimal eigenvalue (resp., singular value) of $M_1\otimes M_2$.


\item $\kappa_\infty(M_1\otimes M_2)\leq \kappa_\infty(A_1\otimes A_2)$.
\end{enumerate}
\end{thm}

\begin{proof}
\begin{enumerate}[i)]
\item By Corollary 1 of \cite{DP}, $|A_1^{-1}|$ dominates $|M_1^{-1}|$ and $|A_2^{-1}|$ dominates $|M_2^{-1}|$.
Hence, $|A_1^{-1}|\otimes |A_2^{-1}|$ dominates $|M_1^{-1}|\otimes |M_2^{-1}|$, and,
since $|(A_1\otimes A_2)^{-1}|=|A_1^{-1}\otimes A_2^{-1}|=|A_1^{-1}|\otimes |A_2^{-1}|$ by \eqref{prop.inv}, $|(A_1\otimes A_2)^{-1}|$ dominates $(M_1\otimes M_2)^{-1}$.

\item Let $B_1$ be an $n\times n$ matrix  and $B_2$ an $m\times m$ matrix. If $\lambda$ is an eigenvalue of $B_1$ and $\mu$ is an eigenvalue of $B_2$, then $\lambda\mu$ is an eigenvalue of $B_1\otimes B_2$ and every eigenvalue of $B_1\otimes B_2$ arises as such a product of eigenvalues of $B_1$ and $B_2$ (see Theorem 4.2.12 of \cite{topics}). By Corollary 2 of \cite{DP}, we have that $\lambda_{\min}(A_1)\leq \lambda_{\min}(M_1)$ and that $\lambda_{\min}(A_2)\leq \lambda_{\min}(M_2)$. Hence, 
\[\lambda_{\min}(M_1\otimes M_2)=\lambda_{\min}(M_1)\lambda_{\min}(M_2)\geq\lambda_{\min}(A_1)\lambda_{\min}(A_2)=\lambda_{\min}(A_1\otimes A_2).\]

The case of singular values is analogous to that of eigenvalues recalling that
every nonzero singular value of $B_1\otimes B_2$ is the product of a singular value of $B_1$ and a singular value of $B_2$ (see Theorem 4.2.15 of \cite{topics}).


\item First, let us see that the infinity norm of the Kronecker product of two matrices $A=(a_{ij})_{1\leq i,j \leq m}$ and $B=(b_{ij})_{1\leq i,j \leq n}$ satisfies that $|| A \otimes B||_{\infty}=|| A||_{\infty}|| B ||_{\infty}$:
\begin{equation}\label{kro.norm}
 || A \otimes B||_{\infty} =\max_{0\leq i \leq nm-1} \sum_{j=1}^m|a_{t+1,j}|(\sum_{k=1}^n |b_{r+1,k}|) , \  \hbox{where }t=\left[\frac{i}{n}\right], r=i-tn . 
\end{equation}

Denoting $R_t:=\sum_{j=1}^m|a_{tj}|$ and $S_r=\sum_{k=1}^n |b_{rk}|$ we can rewrite \eqref{kro.norm} as
\[ || A \otimes B||_{\infty} =\max_{0\leq tn+r \leq nm-1} R_{t+1}S_{r+1}=\max_{1\leq t\leq m}R_{t}\max_{1\leq r\leq n}S_{r}= || A||_{\infty} ||B||_{\infty}. \]
 Hence, the condition number satisfies by \eqref{prop.inv} that
  \begin{align*}
  \kappa_\infty(B_1\otimes B_2)&=|| B_1 \otimes B_2||_{\infty}|| (B_1 \otimes B_2)^{-1}||_{\infty} \\
  &=|| B_1||_{\infty}|| B_2||_{\infty}|| B_1^{-1} ||_{\infty}||B_2^{-1}||_{\infty}=\kappa_\infty(B_1)\kappa_\infty(B_2).
  \end{align*}
  
By Corollary 2 of \cite{DP}, we have that $\kappa_\infty(M_1)\leq \kappa_\infty(A_1)$ and $\kappa_\infty(M_2)\leq \kappa_\infty(A_2)$. So, we conclude that $\kappa_\infty(M_1\otimes M_2)\leq \kappa_\infty(A_1\otimes A_2)$.
\end{enumerate}

\end{proof}

\section{Numerical tests}

In this section two numerical examples illustrating the theoretical results will be presented. The first example
will be constructed  by performing the tensor product of three different NTP bases $u^n=(u_0^n,\ldots,u_n^n)$ 
of the space $\mathcal{P}_n([0,1])$ of polynomials of degree not greater than $n$, which were used in \cite{DP}. A second example will be presented considering the tensor product of rational bases
$r^n=(r_0^n,\ldots,r_n^n)$
constructed from the three NTP bases considered in the first example with positive weights and the tensor product of
rational monomial bases (the monomial basis is TP in $[0,1]$) also with positive weights. In fact, if $u^n$ is a TP basis, it can be checked that
the rational basis $(r_0^n,\ldots,r_n^n)$, $r_i^n(x)=w_iu_i^n(x)/(\sum_{j=0}^n w_ju_j^n(x))$, with  weights $w_i^n>0$, is NTP. The basis $u^n=(b_0^n,\ldots,b_n^n)$ formed by the Bernstein polynomials of degree $n$ (see Example 6 a) in \cite{DP})
is the normalized B-basis of  $\mathcal{P}_n([0,1])$ and the corresponding rational Bernstein basis 
$r^n_B$ defined by $r_i^n(x)=w_ib_i^n(x)/(\sum_{j=0}^n w_jb_j^n(x))$ with $w_i>0$, $i=0,\ldots,n$, is the normalized
B-basis of its spanned space $\langle r^n_B\rangle$.

We will also consider the Said-Ball basis 
$s^n=(s_0^n,\ldots,s_n^n)$ and the DP basis 
$c^n=(c_0^n,\ldots,c_n^n)$, which are both NTP basis.
The Said-Ball basis (see \cite{GS}) is defined by
\begin{align*}
	s_i^n(x)&={\lfloor n/2\rfloor+i\choose i}x^i(1-x)^{\lfloor n/2\rfloor+1},\quad
		0\leq i\leq \lfloor(n-1)/2 \rfloor, 
\end{align*}
$s_i^n(x)=s_{n-i}^n(1-x)$,
		$\lfloor n/2\rfloor+1\leq i\leq n$, and, if $n$ is even 
\[
	s_{n/2}^n(x)={n\choose n/2}x^{n/2}(1-x)^{n/2},
\]
where $\lfloor m\rfloor$ is the greatest
integer less than or equal to $m$. The DP 
basis is given by $c_0^n(x)=(1-x)^n$, $c_n^n(x)=x^n$, $c_i^n(x)=x(1-x)^{n-i},\ 1\leq i\leq\lfloor n/2\rfloor-1$, $c_i^n(x)=x^i(1-x),\ \lfloor (n+1)/2\rfloor+1\leq i\leq n-1$, 
and, if $n$ is even $c_{\frac{n}{2}}^n(x)=1-x^{\frac{n}{2}+1}-(1-x)^{\frac{n}{2}+1}$, and, if $n$ is odd,
\[
	c_{\frac{n-1}{2}}(x)=x(1-x)^{\frac{n+1}{2}}+\frac{1}{2}\left[1-x^{\frac{n+1}{2}+1}-(1-x)^{\frac{n+1}{2}+1}\right],\ 
	c_{\frac{n+1}{2}}(x)=c_{\frac{n-1}{2}}(1-x).
\]
Let $(t_i^n)_{i=1}^{n+1}$ be the sequence of points given 
by $t_i =i/(n+2)$ for $i=1,\ldots,n+1$. Let us consider the Kronecker products of the collocation matrices of the Bernstein,
Said-Ball and DP bases of $\mathcal{P}_n([0,1])$ for
$n=3,4,5$ at $(t_i^n)_{i=1}^{n+1}$ by itself:
$M^n\otimes M^n$, $B_1^n\otimes B_1^n$ and $B_2^n\otimes B_2^n$, respectively. Then, the computation of the eigenvalues and the singular values of these matrices 
have been carried out with Mathematica using a precision of $100$ digits. We can
see the corresponding minimal eigenvalues and singular values in Table \ref{tab:min}. It can be
observed that the minimal eigenvalue, resp. singular value, of $M^n\otimes M^n$ is higher than the
minimal eigenvalue, resp. singular value, of $B_1^n\otimes B_1^n$ 
and $B_2^n\otimes B_2^n$ as Theorem \ref{thm:main} has stated.


\begin{table}[h!]
\centering\begin{tabular}{|r|cc|cc|cc|}
\hline
$n$ & \multicolumn{2}{|c|}{$M^n\otimes M^n$} & \multicolumn{2}{|c|}{$B_1^n\otimes B_1^n$} & \multicolumn{2}{|c|}{$B_2^n\otimes B_2^n$} \\
       & $\lambda_{min}$ & $\sigma_{min}$ & $\lambda_{min}$ & $\sigma_{min}$ & $\lambda_{min}$ & $\sigma_{min}$ \\
 \hline
$3$ & $2.30e-03$ & $2.19e-03$ & $8.28e-04$ & $8.28e-04$ & $3.23e-04$ & $3.20e-04$ \\
$4$ & $3.43e-04$ & $3.23e-04$ & $2.17e-04$ & $1.97e-04$ & $1.92e-05$ & $1.11e-05$ \\
$5$ & $5.10e-05$ & $4.78e-05$ & $1.04e-05$ & $1.03e-05$ & $3.54e-07$ & $2.77e-07$ \\
\hline
\end{tabular}
\caption{The minimal eigenvalue and singular value of $M^n\otimes M^n$, $B_1^n\otimes B_1^n$ and $B_2^n\otimes B_2^n$}\label{tab:min}
\end{table}
We have also computed $k_{\infty}(M^n\otimes M^n)$,  $k_{\infty}(B_1^n\otimes B_1^n)$ and 
$k_{\infty}(B_2^n\otimes B_2^n)$ for $n = 3,4,5$. Table \ref{tab:cond} shows the results. It can be observed that
$k_{\infty}(M^n\otimes M^n)\leq k_{\infty}(B_i^n\otimes B_i^n)$ for $i=1,2$, as it has been shown in 
Theorem \ref{thm:main}.

\begin{table}[h!]
\centering\begin{tabular}{|r|c|c|c|}
\hline
$n$ & $k_{\infty}(M^n\otimes M^n)$ & $k_{\infty}(B_1^n\otimes B_1^n)$ & $k_{\infty}(B_2^n\otimes B_2^n)$ \\
 \hline
$3$ & $5.1883e+02$ & $1.7361e+03$ & $7.1797e+03$ \\
$4$ & $3.9690e+03$ & $6.5610e+03$ &	$1.6080e+05$ \\	
$5$ & $2.5264e+04$ & $1.3949e+05$ & $6.0028e+06$ \\	
\hline
\end{tabular}
\caption{Infinity condition number  $k_{\infty}$ of $M^n\otimes M^n$, $B_1^n\otimes B_1^n$ and $B_2^n\otimes B_2^n$}\label{tab:cond}
\end{table}

As it has been said before, the rational Said-Ball, DP and monomial bases
with positive weights are NTP.
Taking a sequence of positive weights $(w_i^n)_{i=0}^n$ and taking into account that 
$\sum_{j=0}^n w_j^n b_j^n(x)\in\mathcal{P}_n([0,1])$ and that $s^n$, $c^n$ and $m^n=(1,x,\ldots,x^n)$
are bases of $\mathcal{P}_n([0,1])$, then there exists three sequence of
weights $(\overline{w}_i^n)_{i=0}^n$, $(\tilde{w}_i^n)_{i=0}^n$ and $(\hat{w}_i^n)_{i=0}^n$
satisfying
\begin{equation}\label{eq:conv}
\sum_{j=0}^n w_j^n b_j^n(x)=\sum_{j=0}^n \overline{w}_j^n s_j^n(x)=\sum_{j=0}^n \tilde{w}_j^n b_j^n(x)
	=\sum_{j=0}^n \hat{w}_j^n c_j^n(x),\quad x\in [0,1].
\end{equation}
Sequences of positive weights $(w_i^n)_{i=0}^n$ have been randomly generated for $n=3,4,5$, where 
each $w_i^n$ is an integer in the interval $[1, 1000]$, until we
have obtained a sequence such that there exists positive sequences $(\overline{w}_i^n)_{i=0}^n$,
$(\tilde{w}_i^n)_{i=0}^n$ and $(\hat{w}_i^n)_{i=0}^n$
satisfying \eqref{eq:conv}. Then we have the normalized B-basis $r_B$, and the
NTP rational bases of $\langle r_B\rangle$ corresponding to the Said-Ball basis, the DP 
basis and the monomial basis.  So, in the second example we have considered the Kronecker products of the collocation matrices of the generated rational Bernstein,
Said-Ball, DP and monomial bases for
$n=3,4,5$ at $(t_i^n)_{i=1}^{n+1}$ by itself:
$M_T^n=MR^n\otimes MR^n$, $B_{1,T}^n=BR_1^n\otimes BR_1^n$, $B_{2,T}^n=BR_2^n\otimes BR_2^n$ and $B_{3,T}^n=BR_3^n\otimes BR_3^n$, respectively. Then, the computation of the eigenvalues and the singular values of these matrices 
have been carried out with Mathematica using a precision of $100$ digits. We can
see the corresponding minimal eigenvalues and singular values of $M_T^n$, $B_{1,T}^n$ and $B_{3,T}^n$ in Table \ref{tab:minR}. It can be
observed that the minimal eigenvalue, resp. singular value, of $M_T^n$ is higher than the
minimal eigenvalue, resp. singular value, of $B_{1,T}^n$ 
and $B_{3,T}^n$ as Theorem \ref{thm:main} has proved.
\begin{table}[h!]
\centering\begin{tabular}{|r|cc|cc|cc|}
\hline
$n$ & \multicolumn{2}{|c|}{$M_T^n$} & \multicolumn{2}{|c|}{$B_{1,T}^n$} & \multicolumn{2}{|c|}{$B_{3,T}^n$} \\
       & $\lambda_{min}$ & $\sigma_{min}$ & $\lambda_{min}$ & $\sigma_{min}$
        & $\lambda_{min}$ & $\sigma_{min}$ \\
 \hline
$3$ & $1.95e-03$ & $1.74e-03$ &	$4.06e-04$ & $3.78e-04$  &	$4.39e-06$ & $3.82e-6$ \\
$4$ & $2.57e-04$ & $2.05e-04$ & $1.30e-04$ & $1.09e-04$ & 	$8.86e-08$ & $2.35e-08$ \\
$5$ & $4.75e-05$ & $4.36e-05$ &	$8.83e-06$ & $8.66e-06$ & 	$2.60e-10$ & $1.63e-10$ \\
\hline
\end{tabular}
\caption{The minimal eigenvalue and singular value of $M_T^n$, $B_{1,T}^n$ and  $B_{3,T}^n$}\label{tab:minR}
\end{table}

We have also computed $k_{\infty}(M_T^n)$,  $k_{\infty}(B_{1,T}^n)$,
$k_{\infty}(B_{2,T}^n)$ and $k_{\infty}(B_{3,T}^n)$ for $n = 3,4,5$
with Mathematica. The results can be seen in Table \ref{tab:condR}. It can be observed that
$k_{\infty}(M_T^n)\leq k_{\infty}(B_{i,T}^n)$ for $i=1,2,3$ (see 
Theorem \ref{thm:main}).

\begin{table}[h!]
\centering\begin{tabular}{|r|c|c|c|c|}
\hline
$n$ & $k_{\infty}(M_T^n)$ & $k_{\infty}(B_{1,T}^n)$ & $k_{\infty}(B_{2,T}^n)$
&  $k_{\infty}(B_{3,T}^n)$ \\
\hline
$3$ & $8.1049e+02$ & $5.6308e+03$ & $3.5425e+04$ & $5.8525e+05$ \\
$4$ & $7.1105e+03$ & $1.3484e+04$ & $2.0327e+06$ & $1.3229e+08$ \\
$5$ & $3.1318e+04$ & $1.6543e+05$ & $4.0614e+07$ & $1.7440e+10$ \\
\hline
\end{tabular}
\caption{Infinity condition number $k_{\infty}$ of $M_T^n$, $B_{1,T}^n$,  $B_{2,T}^n$ and  $B_{3,T}^n$}\label{tab:condR}
\end{table}


\end{document}